\documentclass[11pt]{article}

\usepackage{xr}
\usepackage{lscape}

\usepackage{amsmath}
\usepackage{amssymb}
\usepackage{amsthm}
\usepackage{enumitem}
\usepackage{titlesec}
\usepackage{graphicx}
\usepackage{caption}
\usepackage{subcaption}

\usepackage{alphalph}

\usepackage[T1]{fontenc}
\usepackage[utf8]{inputenc}

\usepackage{comment}

\usepackage{indentfirst}
\usepackage[top=1.25in,left=1.25in,right=1.25in]{geometry}

\newlength{\widtwo}
\newlength{\wida}
\newlength{\nums}
\newlength{\numone}
\newlength{\widone}
\newlength{\numa}
\newlength{\numc}
\newlength{\numd}
\newtheorem{thm}{Theorem}

\theoremstyle{definition}

\newtheorem{que}[thm]{Question}

\newtheorem{rem}[thm]{Remark}

\numberwithin{equation}{section}

\author{\Large{Riccardo W. Maffucci\footnote{\texttt{riccardo.maffucci@epfl.ch}.}}}
\title{\Large{\uppercase{\bf On polyhedral graphs and their complements}}}

\date{}

\def\ov{\overline}

\begin{document}
\titleformat{\section}
  {\Large\scshape}{\thesection}{1em}{}
\titleformat{\subsection}
  {\large\scshape}{\thesubsection}{1em}{}
\maketitle


\begin{abstract}
	We find all polyhedral graphs such that their complements are still polyhedral. These turn out to be all self-complementary.
\end{abstract}


\section{Introduction}
\subsection{The problem}
The problem we investigate combines two main ideas. Polyhedral graphs, or simply polyhedra, are $3$-connected, planar graphs. This class of graphs is closely related to $3$-dimensional topology and geometry, and the name comes from the fact that they are $1$-skeletons of polyhedral solids in the sense of geometry (Rademacher-Steinitz's Theorem, see e.g. \cite[Theorem 11.6]{harary}). In what follows, we assume polyhedral solids to be convex, and consider them up to topology, i.e., up to their $1$-skeletons being isomorphic graphs. It may be shown using only graph theory that there are only five regular polyhedral solids, namely the Platonic ones \cite[Theorem 1.38]{hahimo}.

Polyhedral graphs have several nice properties. They are the planar graphs that can be embedded in a sphere in a unique way (an observation due to Whitney, see e.g. \cite[Theorem 11.5]{harary}). Specifically, the dual graph is always unique, and duals of polyhedra (in the sense of both graph theory and geometry) are also polyhedra (e.g., \cite[Chapter 11]{harary}). We also record that all their regions (or `faces') are delimited by cycles (elementary closed walks) \cite[Proposition 4.26]{dieste}.




The other idea comes from a classical problem in graph theory, set by Harary: to find all graphs $G$ such that a certain property is verified by both $G$ and its complement graph $\ov{G}$ \cite[Introduction]{coplan}. Our problem is the following.
\begin{que}
	\label{qu}
	Which pairs of complementary graphs $G,\ov{G}$ are both polyhedral?
\end{que}

\begin{figure}[h!]
	\begin{subfigure}{.32\textwidth}
		\centering	\includegraphics[width=3cm,clip=false]{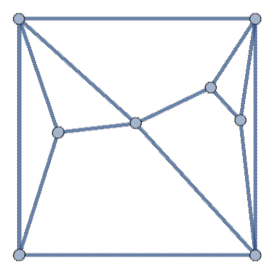}
		\caption{$g_{1408.12}$}
		\label{pic23}
	\end{subfigure}
	\begin{subfigure}{.32\textwidth}
		\centering
		\includegraphics[width=3cm,clip=false]{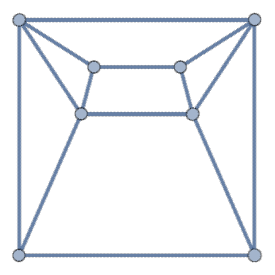}
		\caption{$g_{1408.13}$}
		\label{pic30}
	\end{subfigure}
\begin{subfigure}{.32\textwidth}
	\centering	\includegraphics[width=3cm,clip=false]{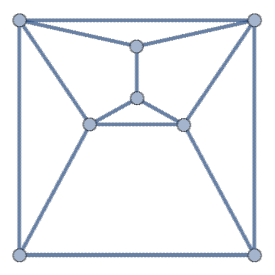}
	\caption{$g_{1408.39}$}
	\label{$g_{31}$}
\end{subfigure}
	\caption{The only three solutions to Question \ref{qu}.}
	\label{fig:01}
\end{figure}

\begin{thm}
\label{thm:1}
There exist exactly three polyhedral graphs such that their complements are polyhedral. These are all self-complementary. They are depicted in Figure \ref{fig:01}.
\end{thm}

All solutions to Question \ref{qu} are $(8,14)$ graphs, of degree sequence
\[4,4,4,4,3,3,3,3.\]
This case promises to be the most interesting, due to the following.
\begin{rem}
	\label{rem:}
	If the polyhedron $G$, its dual, and its complement graphs are all of same order and size, then $G$ is an $(8,14)$ graph. To see this, we impose the following three conditions. If $G$ and $G'$ have the same order, then the number of regions and vertices of $G$ coincide. if $G,\ov{G}$ have the same size, then $q=\frac{1}{2}\binom{p}{2}$. The third condition is Euler's formula for planar graphs. Solving the resulting system, we get $p=8$ and $q=14$.
\end{rem}

\paragraph{Related problems.}
Planar graphs with planar complement were investigated in \cite{bahako}, \cite{coplan}. In \cite[Figure 3.1]{abhill}, we find the only three non-trivial, self-complementary, self-dual graphs. In this figure, graph A is $g_{1408.13}$ of Figure \ref{fig:01}, C is $g_{1408.12}$, while B is not $3$-connected. By the way, $g_{1408.39}$ in Figure \ref{fig:01} is the only self-complementary, non-self-dual polyhedron, as we shall see in section \ref{sec:13}. In other related work, Ando-Kaneko \cite{andkan} investigated the connectivity of complements of $3$-connected graphs.


\subsection{Notation and conventions.}
We say that $G$ is a $(p,q)$ graph when $G$ is a graph of order (number of vertices/points) $p$ and size (number of edges/lines) $q$.
The vertex and edge sets of $G$ are $V(G)$ and $E(G)$ respectively. Vertices are denoted by $v_1,v_2,\dots,v_p$. Their respective degrees are non-negative integers $d_1,d_2,\dots,d_p$. The degree sequence of a graph is indicated by $d_1,d_2,\dots,d_p$, or simply by $d_1d_2\dots d_p$ if there is no possibility of confusion.

A graph of order at least $4$ is $3$-connected if removing any set of $0$, $1$, or $2$ vertices produces a connected graph. The complement of $G$ is denoted by $\ov{G}$, and the dual by $G'$. The letters $\ov{p},\ov{q},p',q'$ indicate their orders and sizes accordingly. The number of regions of a planar graph is $r$, and of its complement and dual (if these are also planar) $\ov{r},r'$. The faces of a polyhedral graph are its regions. The faces are triangular, quadrilateral, pentagonal, ... if they are bounded by a cycle of length $3,4,5,...$ respectively. 



\subsection{Acknowledgements}
The author was supported by Swiss National Science Foundation project 200021\_184927.

\section{Setup}
\label{sec:pre}
It is well-known that if $G$ is a planar graph on at least $9$ vertices, then $\ov{G}$ is non-planar \cite{bahako}. It follows right away that Question \ref{qu} has a finite number of solutions.

On the other hand, we will now show that if $G,\ov{G}$ are polyhedra, then $p\geq 8$. Denote by
\[d_1,d_2,\dots,d_p\]
the (weakly decreasing) degree sequence of $G$. Since the graph is $3$-connected, in particular one has 
\begin{equation}
\label{del}
d_p\geq 3.
\end{equation}
Accordingly, the degree sequence of $\ov{G}$ is
\[p-1-d_p,p-1-d_{p-1},\dots,p-1-d_1,\]
with $p-1-d_1\geq 3$, i.e.
\begin{equation}
\label{Del}
d_1\leq p-4.
\end{equation}
By \eqref{del} and \eqref{Del}, $p\geq 7$. If $p=7$, then $G$ is $3$-regular, impossible due to the handshaking lemma. Therefore, $p\geq 8$, so that ultimately $p=8$.

Further, the above yields $d_1=4$ and $d_8=3$. The handshaking lemma now reduces our cases to those of Table \ref{tab}, assuming w.l.o.g. that $q\leq\ov{q}$. It will thus suffice to inspect the cases $q=12,13,14$. This will be done in the next section.

\begin{table}[h!]
	\centering
$\begin{array}{|c|c|c||c|c|c|}
\hline
\text{deg. sequence of } G&\text{size } q&\text{faces } r&\text{deg. sequence of } \ov{G}&\text{size } \ov{q}&\text{faces } \ov{r}
\\
\hline
33333333&12&6&44444444&16&10
\\
44333333&13&7&44444433&15&9
\\
44443333&14&8&44443333&14&8
\\
\hline
\end{array}$
\caption{Possible degree sequences for solutions of Question \ref{qu}.}
\label{tab}
\end{table}

\section{Completing the proof}
\label{sec:13}
The polyhedra up to $8$ faces are tabulated in \cite{brdu73} and \cite{fede75}. In Appendix \ref{appa}, we collect those up to $14$ edges, for quick reference.

\begin{figure}[h!]
\begin{subfigure}{\nums}
	\centering
	\includegraphics[width=\wida,clip=false]{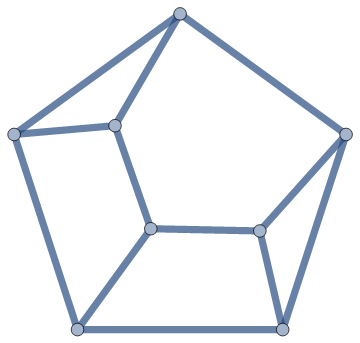}
	\includegraphics[width=\wida,clip=false]{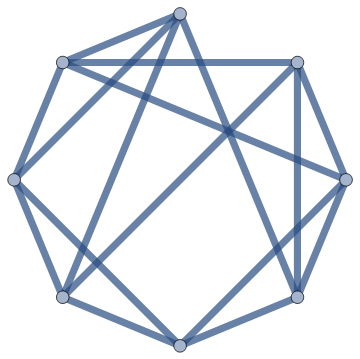}
\end{subfigure}
\begin{subfigure}{\nums}
	\centering
	\includegraphics[width=\wida,clip=false]{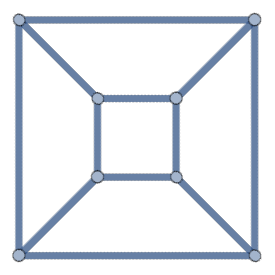}
	\includegraphics[width=\wida,clip=false]{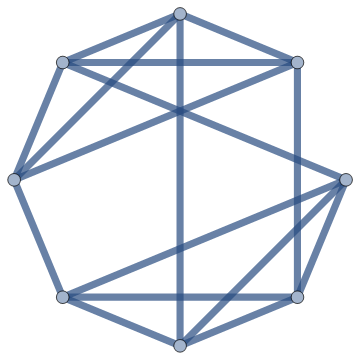}
\end{subfigure}
\caption{The two $(8,12)$ polyhedral graphs, and their complements.}
\label{fig:uno}
\end{figure}

There are only two $(8,12)$ polyhedra (see Figure \ref{fig:12}). For both, the complement is non-planar -- see Figure \ref{fig:uno}.

\begin{figure}[h!]
\begin{subfigure}{\nums}
	\centering
	\includegraphics[width=\wida,clip=false]{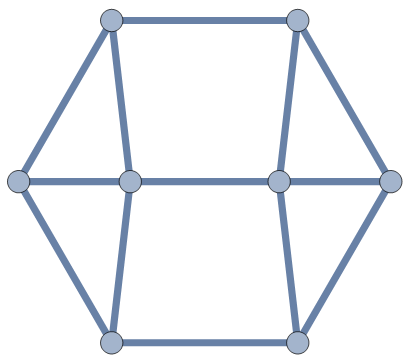}
	\includegraphics[width=\wida,clip=false]{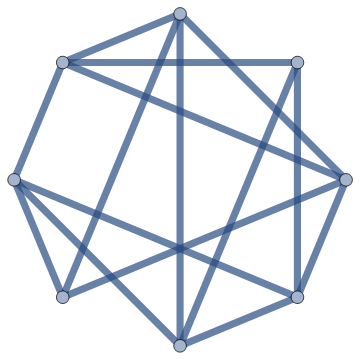}
	\caption{$g_{13.04}$ and its complement}
\end{subfigure}
\begin{subfigure}{\nums}
	\centering
	\includegraphics[width=\wida,clip=false]{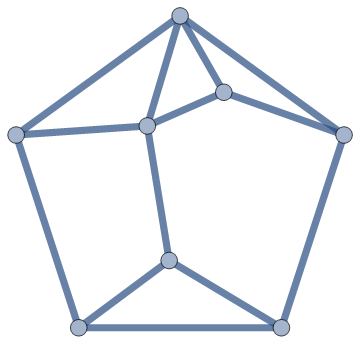}
	\includegraphics[width=\wida,clip=false]{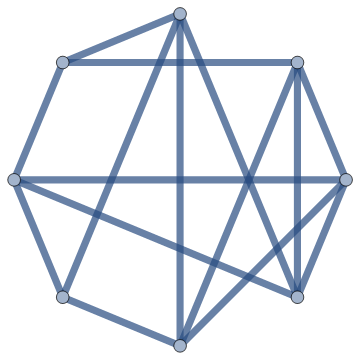}
	\caption{$g_{13.06}$ and its complement}
\end{subfigure}
\begin{subfigure}{\nums}
	\centering
	\includegraphics[width=\wida,clip=false]{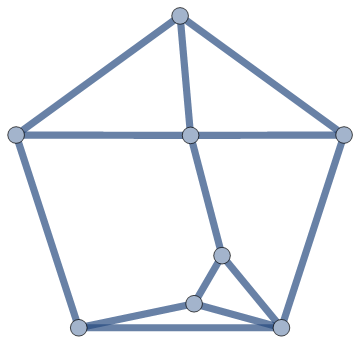}
	\includegraphics[width=\wida,clip=false]{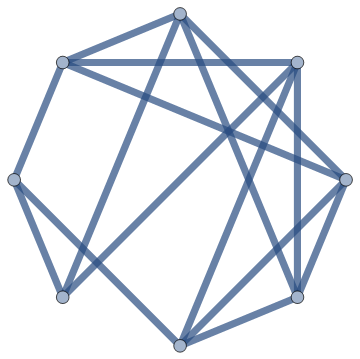}
	\caption{$g_{13.08}$ and its complement}
\end{subfigure}
\begin{subfigure}{\nums}
	\centering
	\includegraphics[width=\wida,clip=false]{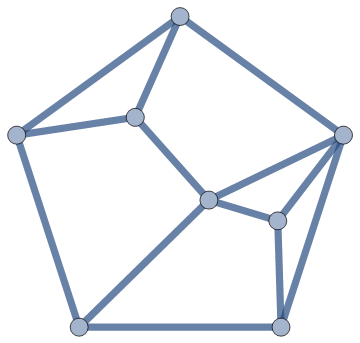}
	\includegraphics[width=\wida,clip=false]{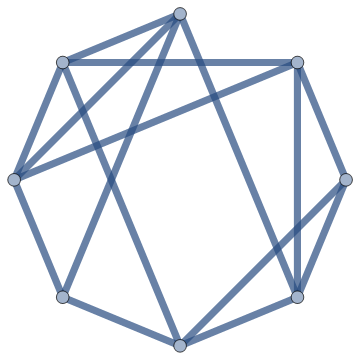}
	\caption{$g_{13.12}$ and its complement}
\end{subfigure}
\begin{subfigure}{\nums}
	\centering
	\includegraphics[width=\wida,clip=false]{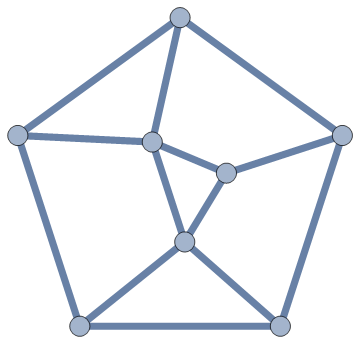}
	\includegraphics[width=\wida,clip=false]{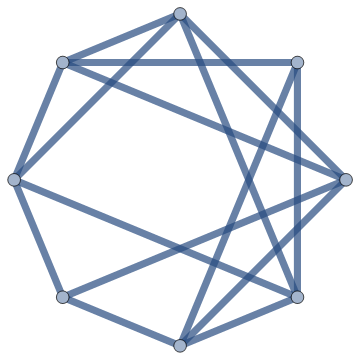}
	\caption{$g_{13.14}$ and its complement}
\end{subfigure}
\begin{subfigure}{\nums}
	\centering
	\includegraphics[width=\wida,clip=false]{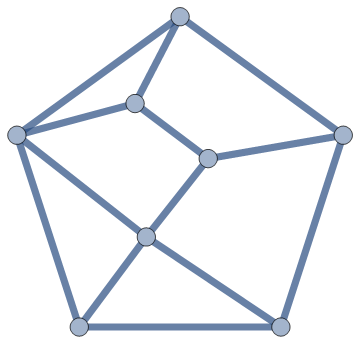}
	\includegraphics[width=\wida,clip=false]{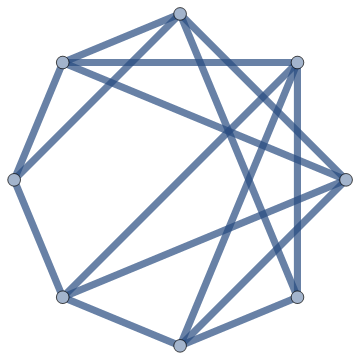}
	\caption{$g_{13.16}$ and its complement}
\end{subfigure}
\begin{subfigure}{\nums}
	\centering
	\includegraphics[width=\wida,clip=false]{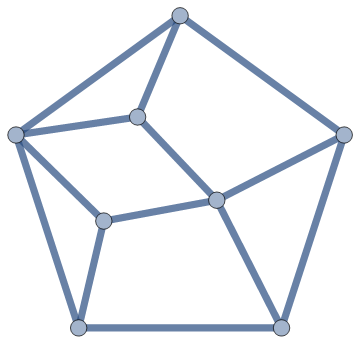}
	\includegraphics[width=\wida,clip=false]{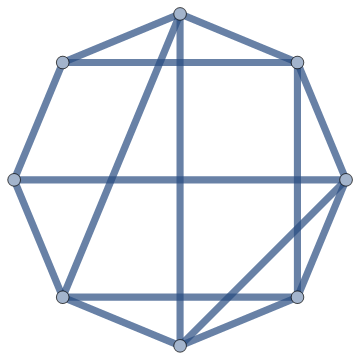}
	\caption{$g_{13.18}$ and its complement}
\end{subfigure}
\begin{subfigure}{\nums}
	\centering
	\includegraphics[width=\wida,clip=false]{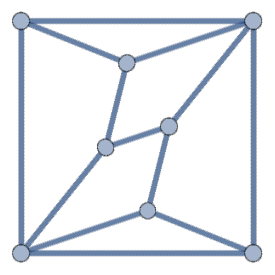}
	\includegraphics[width=\wida,clip=false]{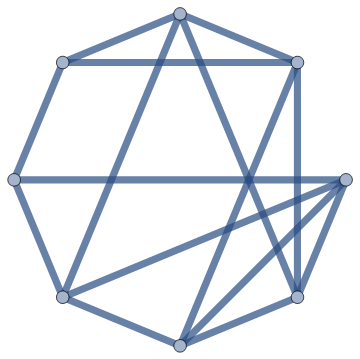}
	\caption{$g_{13.20}$ and its complement}
\end{subfigure}
\begin{subfigure}{\nums}
	\centering
	\includegraphics[width=\wida,clip=false]{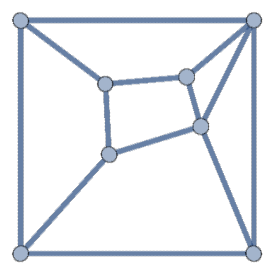}
	\includegraphics[width=\wida,clip=false]{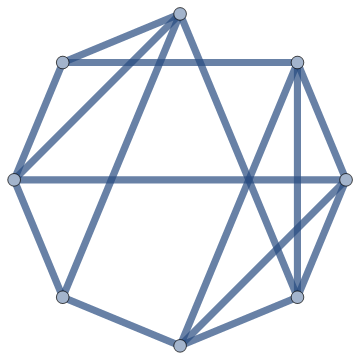}
	\caption{$g_{13.22}$ and its complement}
\end{subfigure}
	\caption{The nine polyhedra of degree sequence $44333333$, and their complement graphs.}
	\label{fig:due}
\end{figure}

There are eleven $(8,13)$ polyhedra, and nine of them have sequence $44333333$ (Figure \ref{fig:13}). In Figure \ref{fig:due}, these are sketched together with their complement graphs. All of the complements are non-planar, hence we discard this case.

\begin{figure}[h!]
	\begin{subfigure}{\nums}
		\centering
		\includegraphics[width=\wida,clip=false]{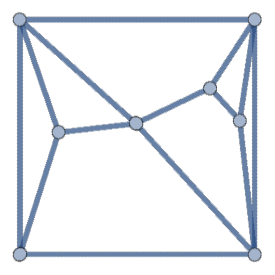}
		\includegraphics[width=\wida,clip=false]{pica23.png}
		\caption{$g_{1408.12}$ and its complement}
	\end{subfigure}
	\begin{subfigure}{\nums}
		\centering
		\includegraphics[width=\wida,clip=false]{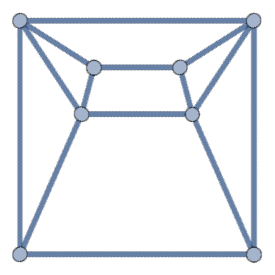}
		\includegraphics[width=\wida,clip=false]{pica30.png}
		\caption{$g_{1408.13}$ and its complement}
	\end{subfigure}
	\begin{subfigure}{\nums}
		\centering
		\includegraphics[width=\wida,clip=false]{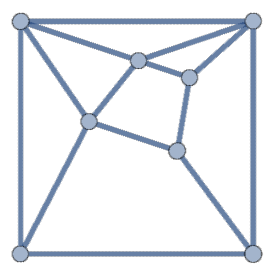}
		\includegraphics[width=\wida,clip=false]{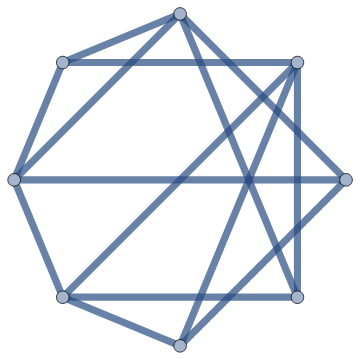}
		\caption{$g_{1408.14}$ and its complement}
	\end{subfigure}
	\begin{subfigure}{\nums}
		\centering
		\includegraphics[width=\wida,clip=false]{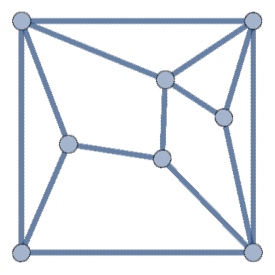}
		\includegraphics[width=\wida,clip=false]{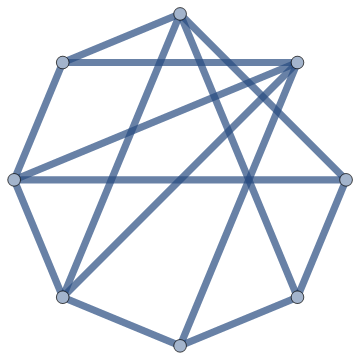}
		\caption{$g_{1408.15}$ and its complement}
	\end{subfigure}
	\begin{subfigure}{\nums}
		\centering
		\includegraphics[width=\wida,clip=false]{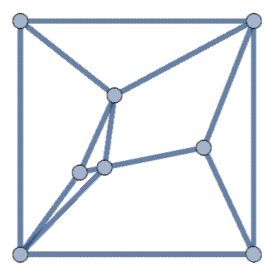}
		\includegraphics[width=\wida,clip=false]{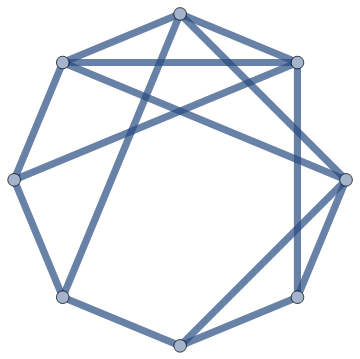}
		\caption{$g_{1408.16}$ and its complement}
	\end{subfigure}
\caption{The self-dual polyhedra of degree sequence $44443333$, and their complement graphs.}
\label{fig:03}
\end{figure}

\begin{figure}[h!]
	\begin{subfigure}{\nums}
		\centering
		\includegraphics[width=\wida,clip=false]{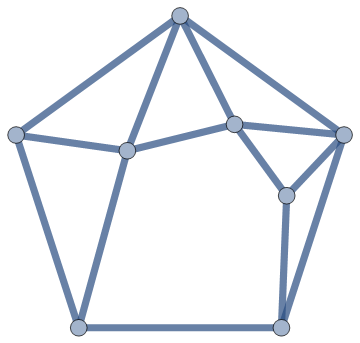}
		\includegraphics[width=\wida,clip=false]{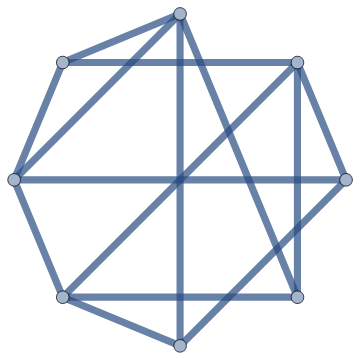}
		\caption{$g_{1408.22}$ and its complement}
	\end{subfigure}
	\begin{subfigure}{\nums}
		\centering
		\includegraphics[width=\wida,clip=false]{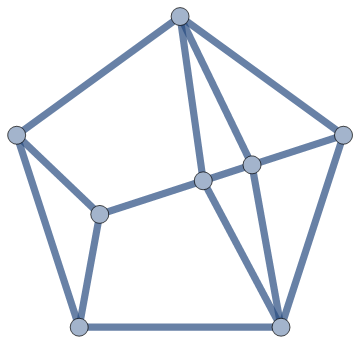}
		\includegraphics[width=\wida,clip=false]{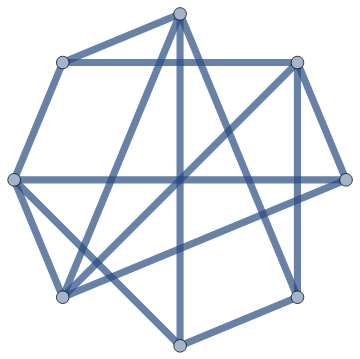}
		\caption{$g_{1408.26}$ and its complement}
	\end{subfigure}
	\begin{subfigure}{\nums}
		\centering
		\includegraphics[width=\wida,clip=false]{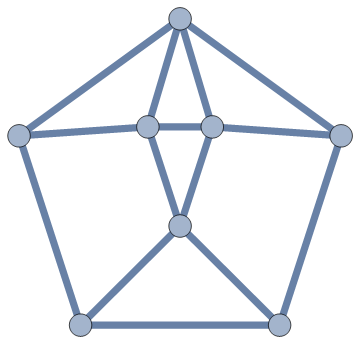}
		\includegraphics[width=\wida,clip=false]{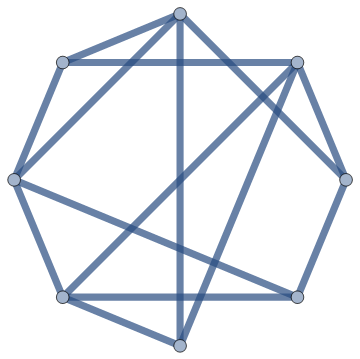}
		\caption{$g_{1408.28}$ and its complement}
	\end{subfigure}
	\begin{subfigure}{\nums}
		\centering
		\includegraphics[width=\wida,clip=false]{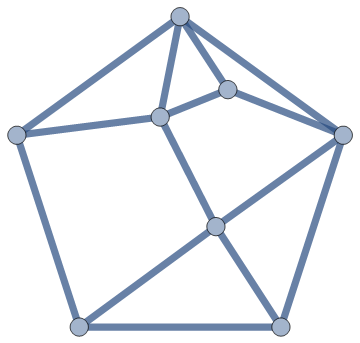}
		\includegraphics[width=\wida,clip=false]{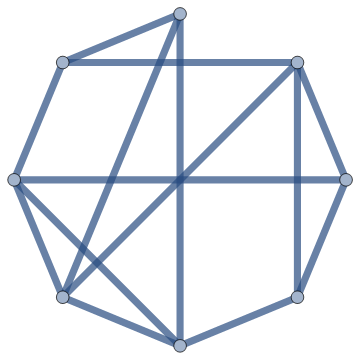}
		\caption{$g_{1408.30}$ and its complement}
	\end{subfigure}
	\begin{subfigure}{\nums}
		\centering
		\includegraphics[width=\wida,clip=false]{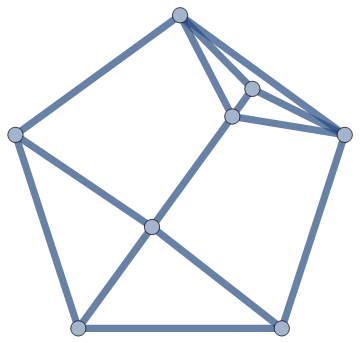}
		\includegraphics[width=\wida,clip=false]{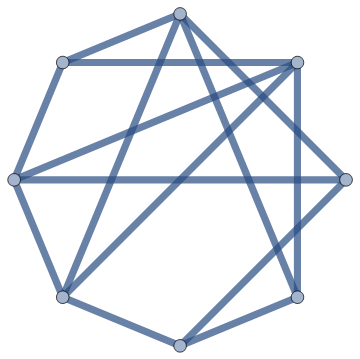}
		\caption{$g_{1408.32}$ and its complement}
	\end{subfigure}
	\begin{subfigure}{\nums}
		\centering
		\includegraphics[width=\wida,clip=false]{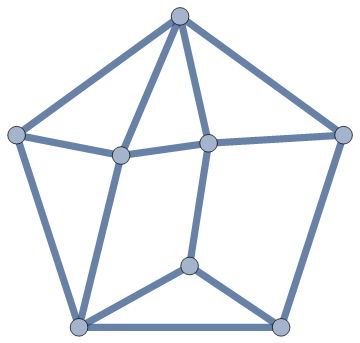}
		\includegraphics[width=\wida,clip=false]{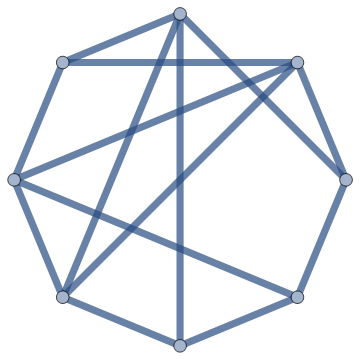}
		\caption{$g_{1408.34}$ and its complement}
	\end{subfigure}
	\begin{subfigure}{\nums}
		\centering
		\includegraphics[width=\wida,clip=false]{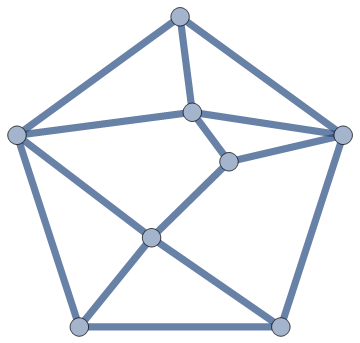}
		\includegraphics[width=\wida,clip=false]{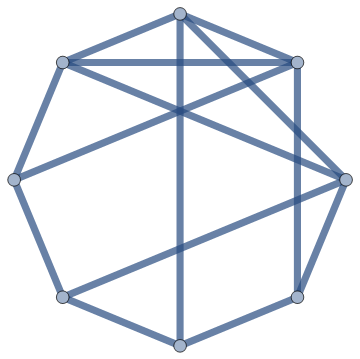}
		\caption{$g_{1408.36}$ and its complement}
	\end{subfigure}
	\begin{subfigure}{\nums}
		\centering
		\includegraphics[width=\wida,clip=false]{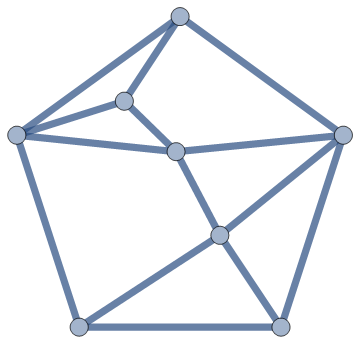}
		\includegraphics[width=\wida,clip=false]{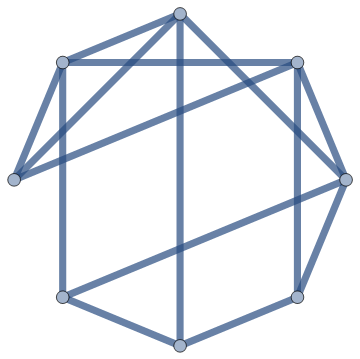}
		\caption{$g_{1408.38}$ and its complement}
	\end{subfigure}
	\begin{subfigure}{\nums}
		\centering
		\includegraphics[width=\wida,clip=false]{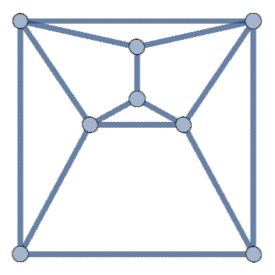}
		\includegraphics[width=\wida,clip=false]{pica31.png}
		\caption{$g_{1408.39}$ and its complement}
	\end{subfigure}
	\begin{subfigure}{\nums}
		\centering
		\includegraphics[width=\wida,clip=false]{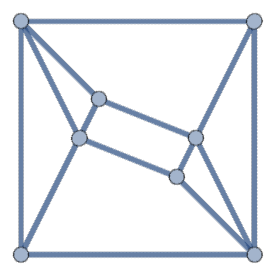}
		\includegraphics[width=\wida,clip=false]{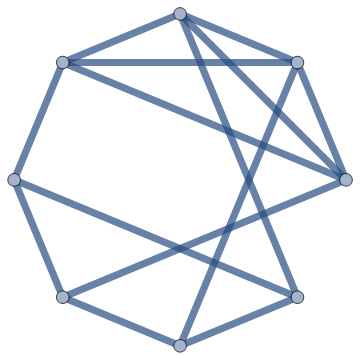}
		\caption{$g_{1408.40}$ and its complement}
	\end{subfigure}
	\begin{subfigure}{\nums}
		\centering
		\includegraphics[width=\wida,clip=false]{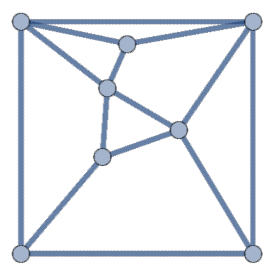}
		\includegraphics[width=\wida,clip=false]{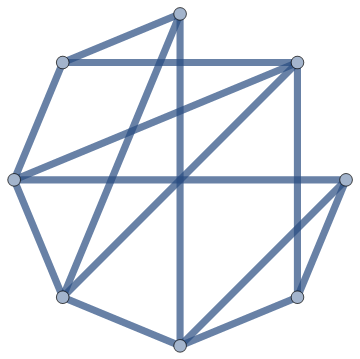}
		\caption{$g_{1408.41}$ and its complement}
	\end{subfigure}
	\begin{subfigure}{\nums}
		\centering
		\includegraphics[width=\wida,clip=false]{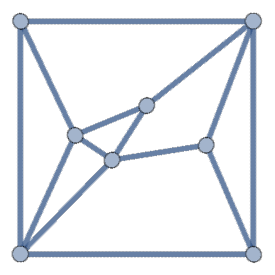}
		\includegraphics[width=\wida,clip=false]{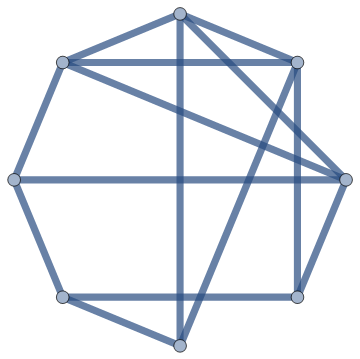}
		\caption{$g_{1408.42}$ and its complement}
	\end{subfigure}
\caption{The non-self-dual polyhedra of degree sequence $44443333$, and their complement graphs.}
\label{fig:04}
\end{figure}

The forty-two $(8,14)$ polyhedra may be found in Figures \ref{fig:14.2} and \ref{fig:14.3}. Exactly seventeen of them have degree sequence $44443333$. These are collected in Figures \ref{fig:03} (self-duals) and \ref{fig:04} (non-self-duals). We find three solutions to Question \ref{qu}, namely graphs $g_{1408.12}$, $g_{1408.13}$, and $g_{1408.39}$.

\begin{figure}[h!]
	\begin{subfigure}{0.48\textwidth}
		\centering
		\includegraphics[width=4.5cm,clip=false]{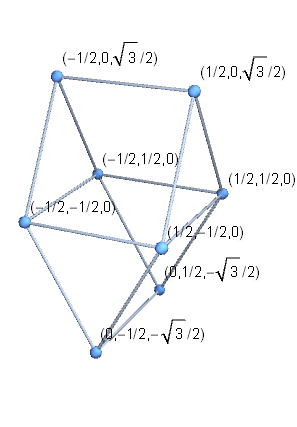}
	\end{subfigure}
	\begin{subfigure}{0.48\textwidth}
		\centering
		\includegraphics[width=4.5cm,clip=false]{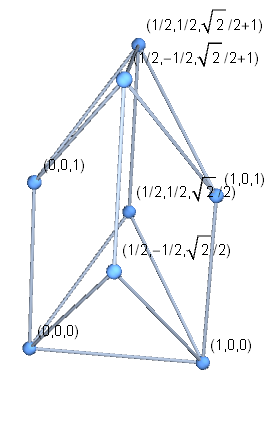}
	\end{subfigure}
	\caption{Embeddings of $g_{1408.39}$ and its dual in $3$-dimensional space, with all edges of unit length.}
	\label{fig:ul}
\end{figure}

As these lists of polyhedra are all exhaustive \cite{brdu73,fede75}, these three are the only solutions to Question \ref{qu}. The proof of Theorem \ref{thm:1} is now complete. We note that $g_{1408.39}$ is the only self-complementary, non-self-dual polyhedron. As a remark of a different flavour, $g_{1408.39}$ and its dual may be embedded in $3$-dimensional space so that all $14$ edges have unit length, as in Figure \ref{fig:ul}.

\clearpage
\appendix
\section{Tables of polyhedra}
\label{appa}
We choose the following ordering. Firstly, the tables are according to increasing size, rather than order. This is due to two main, related reasons. Understanding $(p,q)$ polyhedra of $r=2+q-p$ faces (Euler's formula), $p>r$, is no harder than studying the $(r,q)$ of $p$ faces, and then passing to the duals. In this sense, the complexity grows with $q$ rather than $p$. Moreover, in this way each table lists dual pairs of polyhedra together (as they have the same size, but not necessarily the same order).

Our next criteria for ordering, size being equal, is by increasing number of vertices. We have the inequalities
\begin{equation*}
\frac{q+6}{3}\leq p\leq\frac{2q}{3}
\end{equation*}
due to planarity ($q\leq 3p-6$), and $3$-connectivity (implying $\delta(G)\geq 3$, hence $q\geq 3p/2$ via the handshaking lemma).

Self-duals are listed before dual pairs. All the above being equal, we sort by decreasing highest degrees of vertices, and then by decreasing highest degrees of vertices of the dual. If all the said criteria are equal, we list arbitrarily.


\subsection{Size $q\leq 12$}
\begin{figure}[h!]
	\begin{subfigure}{\numone}
		\centering
		\includegraphics[width=\widone,clip=false]{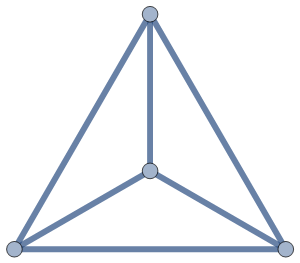}
		\caption{$g_{06}$}
	\end{subfigure}
	\begin{subfigure}{\numone}
		\centering
		\includegraphics[width=\widone,clip=false]{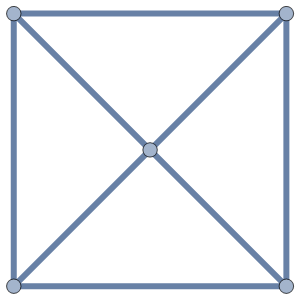}
		\caption{$g_{08}$}
	\end{subfigure}
	\begin{subfigure}{\numone}
		\centering
		\includegraphics[width=\widone,clip=false]{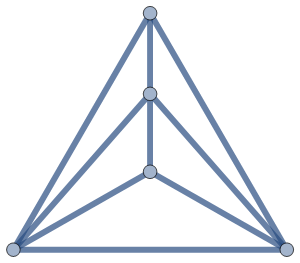}
		\includegraphics[width=\widone,clip=false]{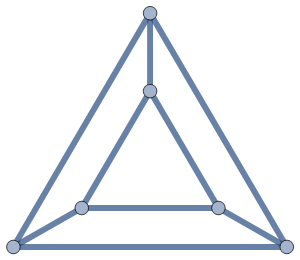}
		\caption{$g_{09.01}$ and its dual}
	\end{subfigure}
	\begin{subfigure}{\numone}
		\centering
		\includegraphics[width=\widone,clip=false]{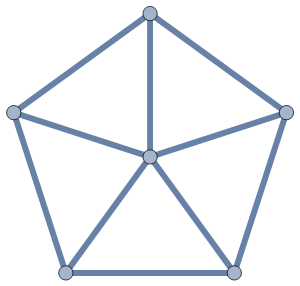}
		\caption{$g_{10.01}$}
	\end{subfigure}
	\begin{subfigure}{\numone}
		\centering
		\includegraphics[width=\widone,clip=false]{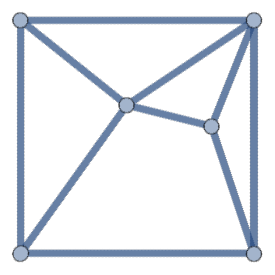}
		\caption{$g_{10.02}$}
	\end{subfigure}
	\caption{The $6$ polyhedra with $q\leq 10$.}
	\label{fig:10}
\end{figure}

\begin{figure}[h!]
	\begin{subfigure}{0.49\textwidth}
		\centering
		\includegraphics[width=\widtwo,clip=false]{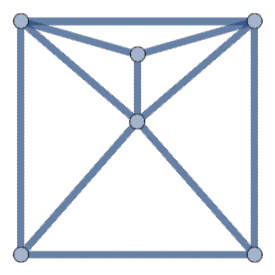}
		\includegraphics[width=\widtwo,clip=false]{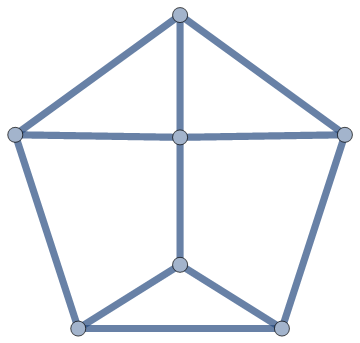}
		\caption{$g_{11.01}$ and its dual}
	\end{subfigure}
	\begin{subfigure}{0.49\textwidth}
		\centering
		\includegraphics[width=\widtwo,clip=false]{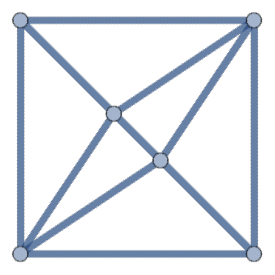}
		\includegraphics[width=\widtwo,clip=false]{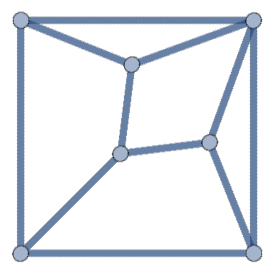}
		\caption{$g_{11.03}$ and its dual}
	\end{subfigure}
	\caption{The $4$ polyhedra with $q=11$.}
	\label{fig:11}
\end{figure}

\begin{figure}[h!]
	\begin{subfigure}{\numc}
		\centering
		\includegraphics[width=\widtwo,clip=false]{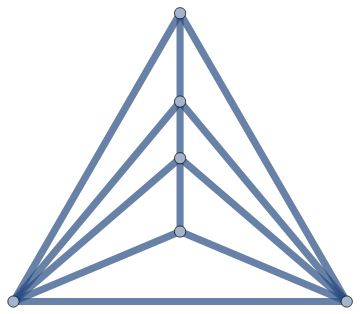}
		\includegraphics[width=\widtwo,clip=false]{h120601.png}
		\caption{$g_{1206.01}$ and its dual}
	\end{subfigure}
	\begin{subfigure}{\numc}
		\centering
		\includegraphics[width=\widtwo,clip=false]{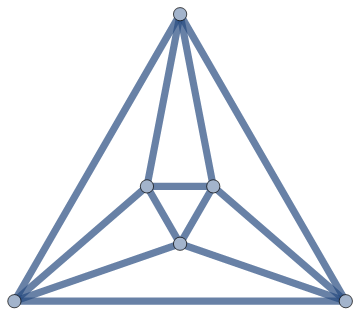}
		\includegraphics[width=\widtwo,clip=false]{h120602.png}
		\caption{$g_{1206.02}$ and its dual}
	\end{subfigure}
	\begin{subfigure}{\numa}
		\centering
		\includegraphics[width=\wida,clip=false]{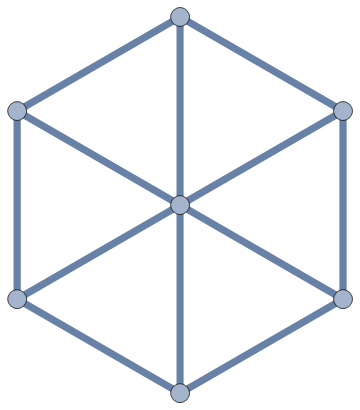}
		\caption{$g_{1207.01}$}
	\end{subfigure}
	\begin{subfigure}{\numa}
		\centering
		\includegraphics[width=\wida,clip=false]{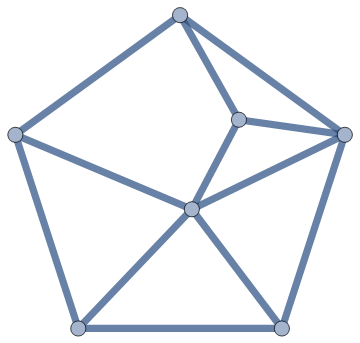}
		\caption{$g_{1207.02}$}
	\end{subfigure}
	\begin{subfigure}{\numa}
		\centering
		\includegraphics[width=\wida,clip=false]{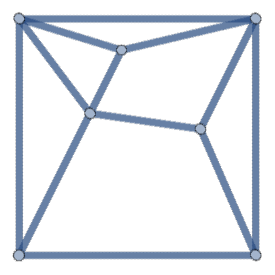}
		\caption{$g_{1207.03}$}
	\end{subfigure}
	\begin{subfigure}{\numa}
		\centering
		\includegraphics[width=\wida,clip=false]{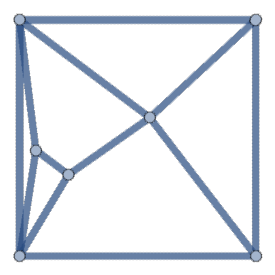}
		\caption{$g_{1207.04}$}
	\end{subfigure}
	\begin{subfigure}{\numa}
		\centering
		\includegraphics[width=\wida,clip=false]{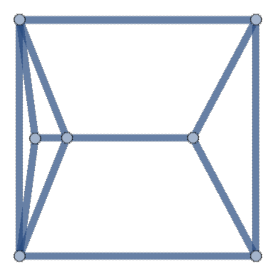}
		\caption{$g_{1207.05}$}
	\end{subfigure}
	\begin{subfigure}{\numa}
		\centering
		\includegraphics[width=\wida,clip=false]{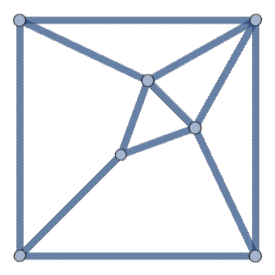}
		\caption{$g_{1207.06}$}
	\end{subfigure}
	\begin{subfigure}{\numc}
		\centering
		\includegraphics[width=\wida,clip=false]{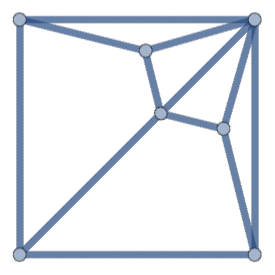}
		\includegraphics[width=\wida,clip=false]{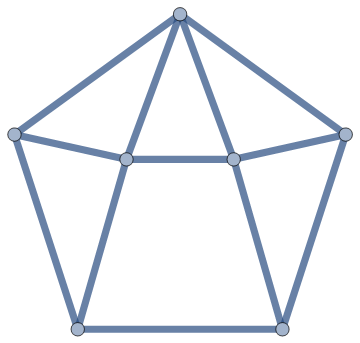}
		\caption{$g_{1207.07}$ and its dual}
	\end{subfigure}
	\caption{The $12$ polyhedra with $q=12$.}
	\label{fig:12}
\end{figure}

\clearpage
\subsection{Size $q=13$}
\begin{figure}[h!]
	\begin{subfigure}{\numd}
		\centering
		\includegraphics[width=\wida,clip=false]{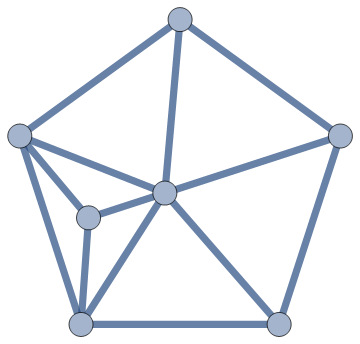}
		\includegraphics[width=\wida,clip=false]{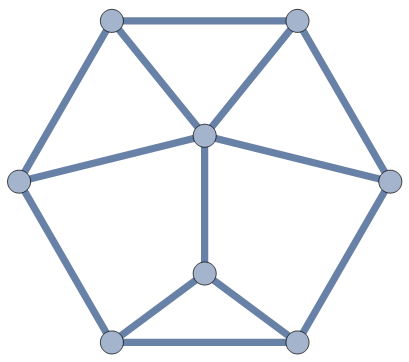}
		\caption{$g_{13.01}$ and its dual}
	\end{subfigure}
	\begin{subfigure}{\numd}
		\centering
		\includegraphics[width=\wida,clip=false]{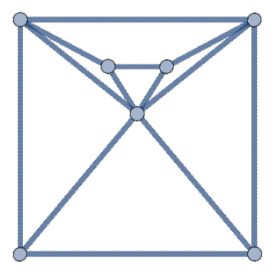}
		\includegraphics[width=\wida,clip=false]{qa5.png}
		\caption{$g_{13.03}$ and its dual}
	\end{subfigure}
	\begin{subfigure}{\numd}
		\centering
		\includegraphics[width=\wida,clip=false]{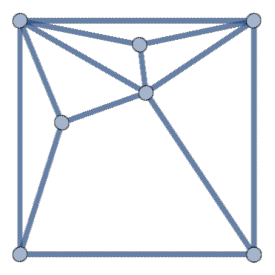}
		\includegraphics[width=\wida,clip=false]{qb121.png}
		\caption{$g_{13.05}$ and its dual}
	\end{subfigure}
	\begin{subfigure}{\numd}
		\centering
		\includegraphics[width=\wida,clip=false]{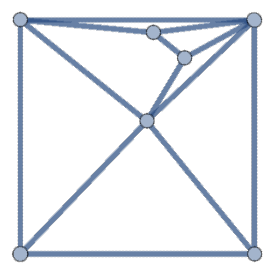}
		\includegraphics[width=\wida,clip=false]{qb122.png}
		\caption{$g_{13.07}$ and its dual}
	\end{subfigure}
	\begin{subfigure}{\numd}
		\centering
		\includegraphics[width=\wida,clip=false]{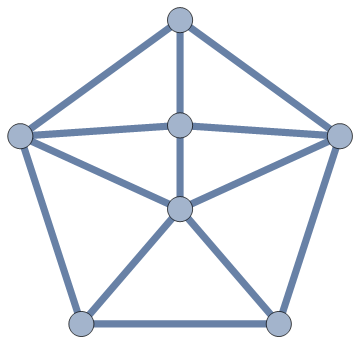}
		\includegraphics[width=\wida,clip=false]{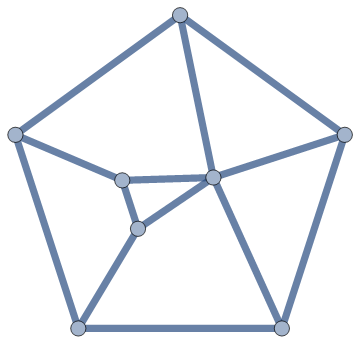}
		\caption{$g_{13.09}$ and its dual}
	\end{subfigure}
	\begin{subfigure}{\numd}
		\centering
		\includegraphics[width=\wida,clip=false]{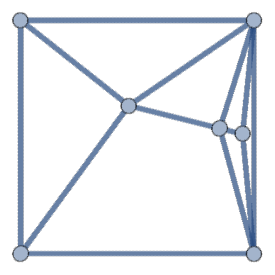}
		\includegraphics[width=\wida,clip=false]{qb5.png}
		\caption{$g_{13.11}$ and its dual}
	\end{subfigure}
	\begin{subfigure}{\numd}
		\centering
		\includegraphics[width=\wida,clip=false]{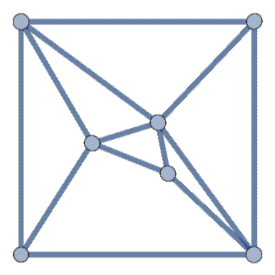}
		\includegraphics[width=\wida,clip=false]{qb61.png}
		\caption{$g_{13.13}$ and its dual}
	\end{subfigure}
	\begin{subfigure}{\numd}
		\centering
		\includegraphics[width=\wida,clip=false]{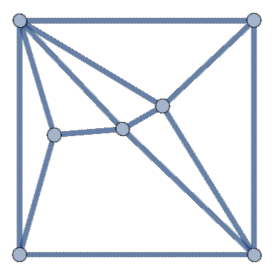}
		\includegraphics[width=\wida,clip=false]{qb123.png}
		\caption{$g_{13.15}$ and its dual}
	\end{subfigure}
	\begin{subfigure}{\numd}
		\centering
		\includegraphics[width=\wida,clip=false]{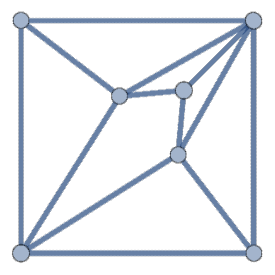}
		\includegraphics[width=\wida,clip=false]{qb161.png}
		\caption{$g_{13.17}$ and its dual}
	\end{subfigure}
	\begin{subfigure}{\numd}
		\centering
		\includegraphics[width=\wida,clip=false]{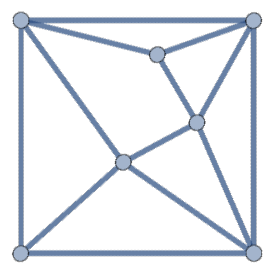}
		\includegraphics[width=\wida,clip=false]{qc2.png}
		\caption{$g_{13.19}$ and its dual}
	\end{subfigure}
	\begin{subfigure}{\numd}
		\centering
		\includegraphics[width=\wida,clip=false]{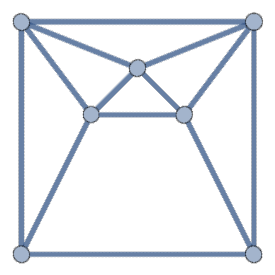}
		\includegraphics[width=\wida,clip=false]{qc5.png}
		\caption{$g_{13.21}$ and its dual}
	\end{subfigure}
	\caption{The $22$ polyhedra with $q=13$.}
	\label{fig:13}
\end{figure}

\clearpage
\subsection{Size $q=14$}
\begin{figure}[h!]
\begin{subfigure}{\numd}
\centering
\includegraphics[width=\wida,clip=false]{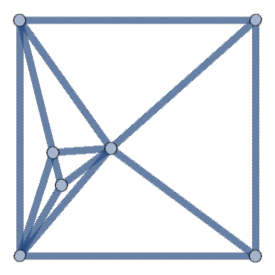}
\includegraphics[width=\wida,clip=false]{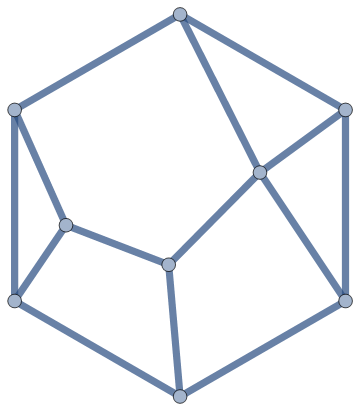}
\caption{$g_{1407.01}$ and its dual}
\end{subfigure}
\begin{subfigure}{\numd}
\centering
\includegraphics[width=\wida,clip=false]{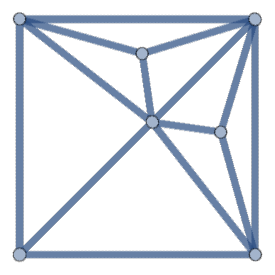}
\includegraphics[width=\wida,clip=false]{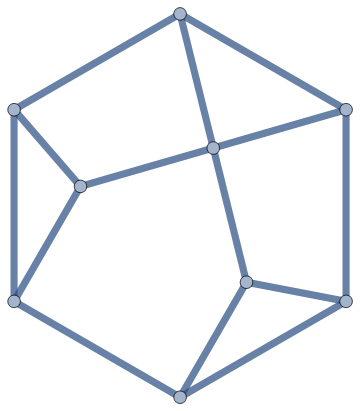}
\caption{$g_{1407.03}$ and its dual}
\end{subfigure}
\begin{subfigure}{\numd}
\centering
\includegraphics[width=\wida,clip=false]{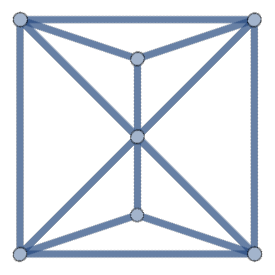}
\includegraphics[width=\wida,clip=false]{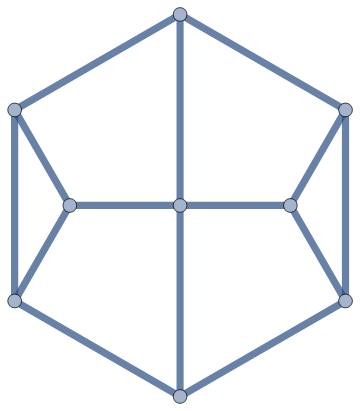}
\caption{$g_{1407.05}$ and its dual}
\end{subfigure}
\begin{subfigure}{\numd}
\centering
\includegraphics[width=\wida,clip=false]{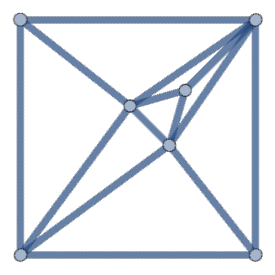}
\includegraphics[width=\wida,clip=false]{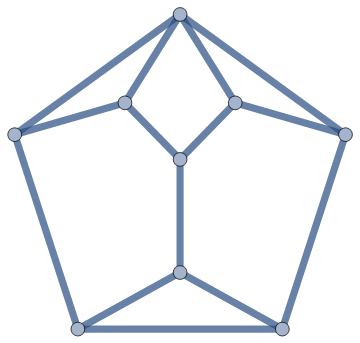}
\caption{$g_{1407.07}$ and its dual}
\end{subfigure}
\begin{subfigure}{\numd}
\centering
\includegraphics[width=\wida,clip=false]{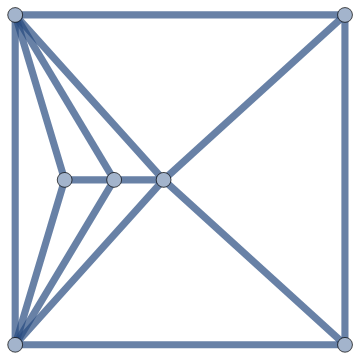}
\includegraphics[width=\wida,clip=false]{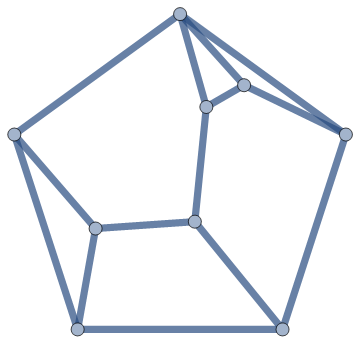}
\caption{$g_{1407.09}$ and its dual}
\end{subfigure}
\begin{subfigure}{\numd}
\centering
\includegraphics[width=\wida,clip=false]{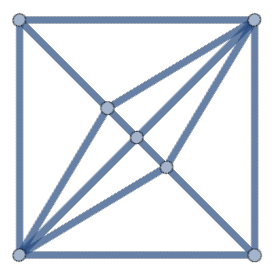}
\includegraphics[width=\wida,clip=false]{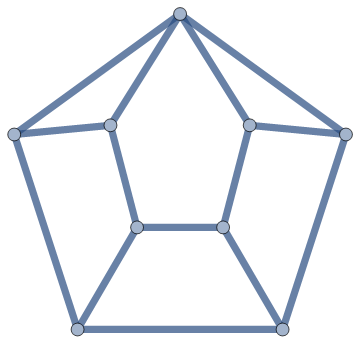}
\caption{$g_{1407.11}$ and its dual}
\end{subfigure}
\begin{subfigure}{\numd}
\centering
\includegraphics[width=\wida,clip=false]{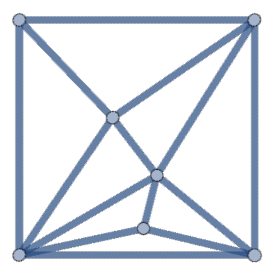}
\includegraphics[width=\wida,clip=false]{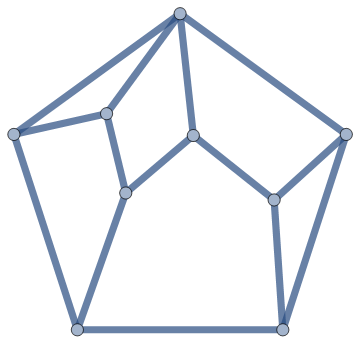}
\caption{$g_{1407.13}$ and its dual}
\end{subfigure}
\begin{subfigure}{\numd}
\centering
\includegraphics[width=\wida,clip=false]{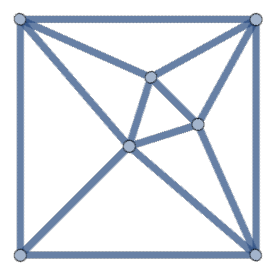}
\includegraphics[width=\wida,clip=false]{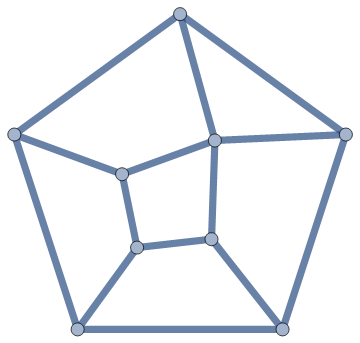}
\caption{$g_{1407.15}$ and its dual}
\end{subfigure}
\caption{The $16$ polyhedra with $q=14$ and $p=7,9$.}
\label{fig:14.1}
\end{figure}

\begin{figure}[h!]
	\begin{subfigure}{\numa}
		\centering
		\includegraphics[width=\wida,clip=false]{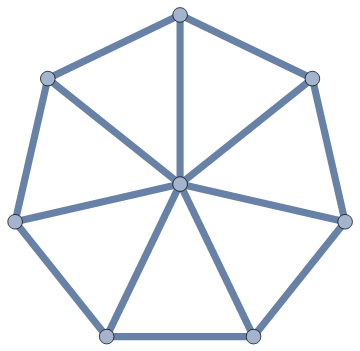}
		\caption{$g_{1408.01}$}
	\end{subfigure}
	\begin{subfigure}{\numa}
		\centering
		\includegraphics[width=\wida,clip=false]{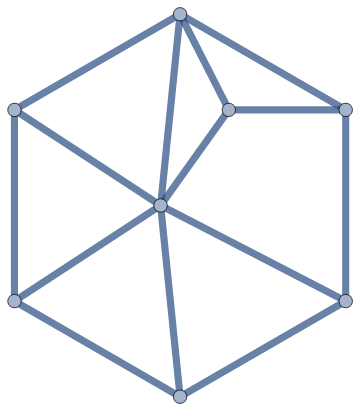}
		\caption{$g_{1408.02}$}
	\end{subfigure}
	\begin{subfigure}{\numa}
		\centering
		\includegraphics[width=\wida,clip=false]{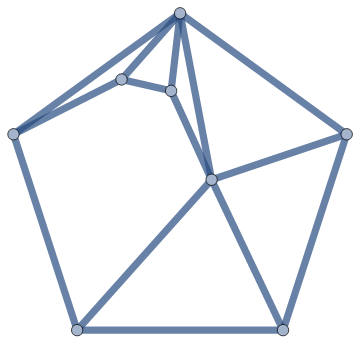}
		\caption{$g_{1408.03}$}
	\end{subfigure}
	\begin{subfigure}{\numa}
		\centering
		\includegraphics[width=\wida,clip=false]{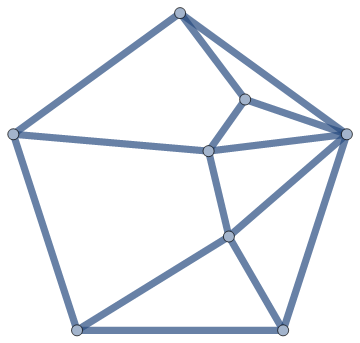}
		\caption{$g_{1408.04}$}
	\end{subfigure}
	\begin{subfigure}{\numa}
		\centering
		\includegraphics[width=\wida,clip=false]{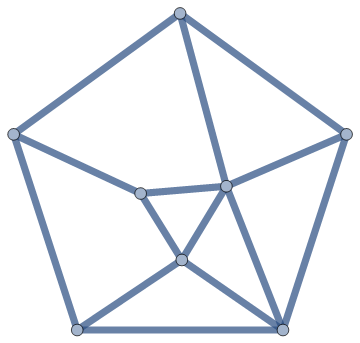}
		\caption{$g_{1408.05}$}
	\end{subfigure}
	\begin{subfigure}{\numa}
		\centering
		\includegraphics[width=\wida,clip=false]{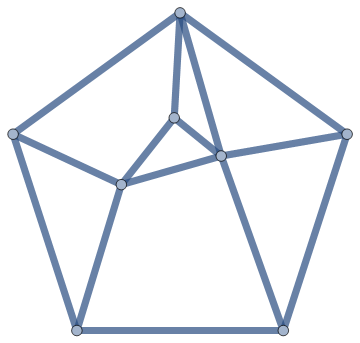}
		\caption{$g_{1408.06}$}
	\end{subfigure}
	\begin{subfigure}{\numa}
		\centering
		\includegraphics[width=\wida,clip=false]{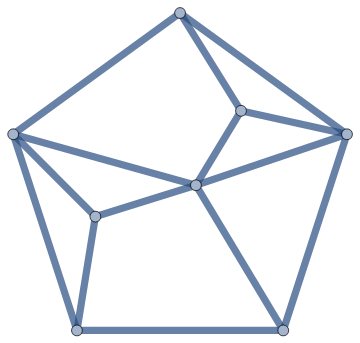}
		\caption{$g_{1408.07}$}
	\end{subfigure}
	\begin{subfigure}{\numa}
		\centering
		\includegraphics[width=\wida,clip=false]{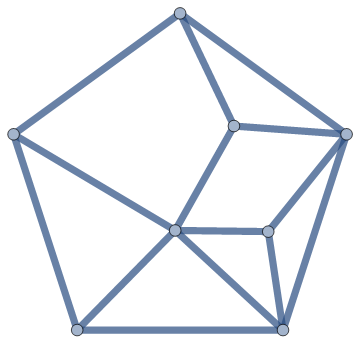}
		\caption{$g_{1408.08}$}
	\end{subfigure}
	\begin{subfigure}{\numa}
		\centering
		\includegraphics[width=\wida,clip=false]{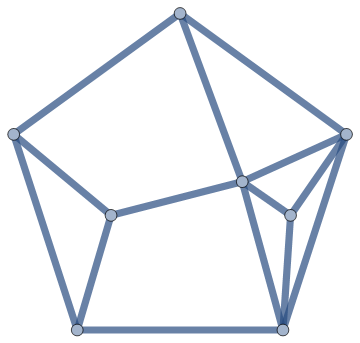}
		\caption{$g_{1408.09}$}
	\end{subfigure}
	\begin{subfigure}{\numa}
		\centering
		\includegraphics[width=\wida,clip=false]{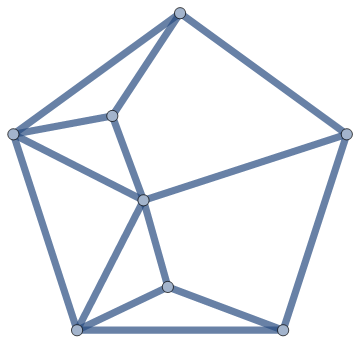}
		\caption{$g_{1408.10}$}
	\end{subfigure}
	\begin{subfigure}{\numa}
		\centering
		\includegraphics[width=\wida,clip=false]{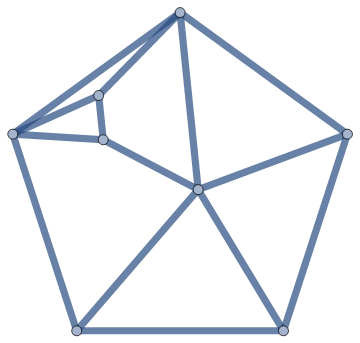}
		\caption{$g_{1408.11}$}
	\end{subfigure}
	\begin{subfigure}{\numa}
		\centering
		\includegraphics[width=\wida,clip=false]{pica23.png}
		\caption{$g_{1408.12}$}
	\end{subfigure}
	\begin{subfigure}{\numa}
		\centering
		\includegraphics[width=\wida,clip=false]{pica30.png}
		\caption{$g_{1408.13}$}
	\end{subfigure}
	\begin{subfigure}{\numa}
		\centering
		\includegraphics[width=\wida,clip=false]{pic15.png}
		\caption{$g_{1408.14}$}
	\end{subfigure}
	\begin{subfigure}{\numa}
		\centering
		\includegraphics[width=\wida,clip=false]{pic38.png}
		\caption{$g_{1408.15}$}
	\end{subfigure}
	\begin{subfigure}{\numa}
		\centering
		\includegraphics[width=\wida,clip=false]{pic12.png}
		\caption{$g_{1408.16}$}
	\end{subfigure}
	\caption{The $16$ self-dual polyhedra with $q=14$ and $p=8$.}
	\label{fig:14.2}
\end{figure}

\begin{figure}[h!]
	\begin{subfigure}{\numa}
		\centering
		\includegraphics[width=\wida,clip=false]{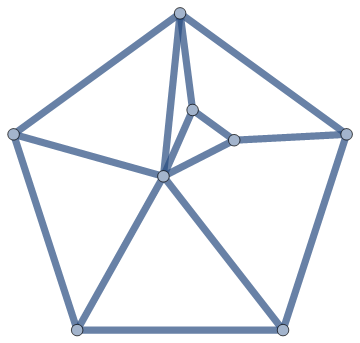}
		\includegraphics[width=\wida,clip=false]{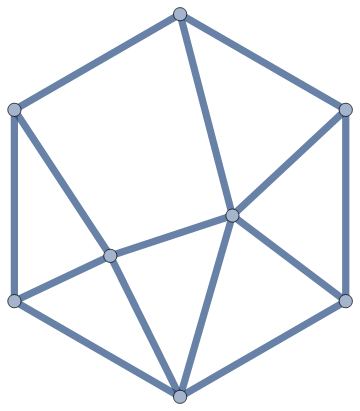}
		\caption{$g_{1408.17}$ and its dual}
	\end{subfigure}
	\begin{subfigure}{\numa}
		\centering
		\includegraphics[width=\wida,clip=false]{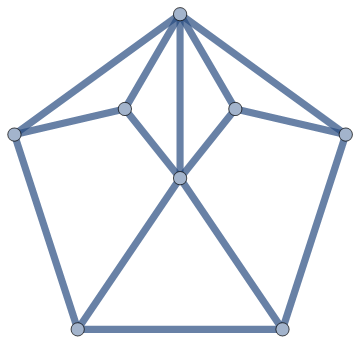}
		\includegraphics[width=\wida,clip=false]{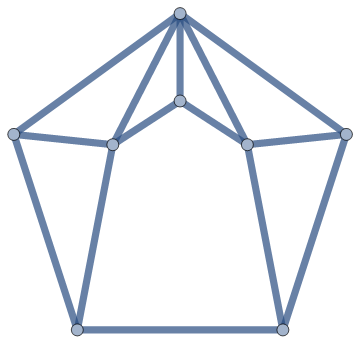}
		\caption{$g_{1408.19}$ and its dual}
	\end{subfigure}
	\begin{subfigure}{\numa}
		\centering
		\includegraphics[width=\wida,clip=false]{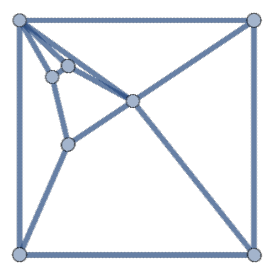}
		\includegraphics[width=\wida,clip=false]{pic17.png}
		\caption{$g_{1408.21}$ and its dual}
	\end{subfigure}
	\begin{subfigure}{\numa}
		\centering
		\includegraphics[width=\wida,clip=false]{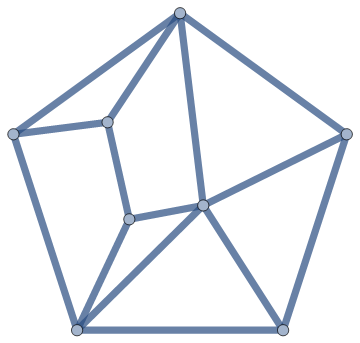}
		\includegraphics[width=\wida,clip=false]{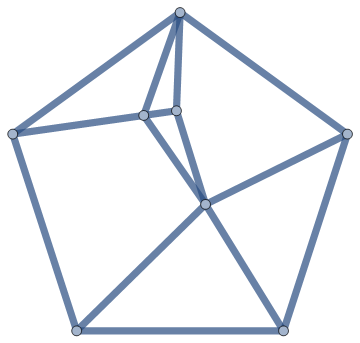}
		\caption{$g_{1408.23}$ and its dual}
	\end{subfigure}
	\begin{subfigure}{\numa}
		\centering
		\includegraphics[width=\wida,clip=false]{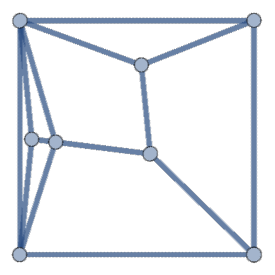}
		\includegraphics[width=\wida,clip=false]{pic03.png}
		\caption{$g_{1408.25}$ and its dual}
	\end{subfigure}
	\begin{subfigure}{\numa}
		\centering
		\includegraphics[width=\wida,clip=false]{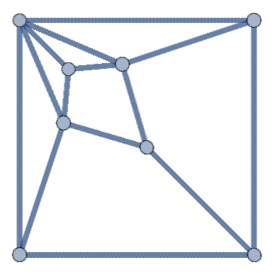}
		\includegraphics[width=\wida,clip=false]{pic13.png}
		\caption{$g_{1408.27}$ and its dual}
	\end{subfigure}
	\begin{subfigure}{\numa}
		\centering
		\includegraphics[width=\wida,clip=false]{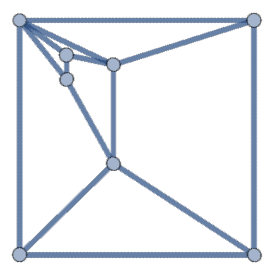}
		\includegraphics[width=\wida,clip=false]{pic19.png}
		\caption{$g_{1408.29}$ and its dual}
	\end{subfigure}
	\begin{subfigure}{\numa}
		\centering
		\includegraphics[width=\wida,clip=false]{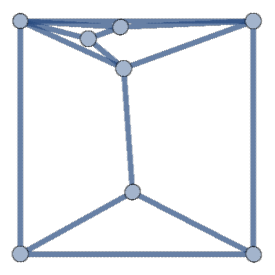}
		\includegraphics[width=\wida,clip=false]{pic33.png}
		\caption{$g_{1408.31}$ and its dual}
	\end{subfigure}
	\begin{subfigure}{\numa}
		\centering
		\includegraphics[width=\wida,clip=false]{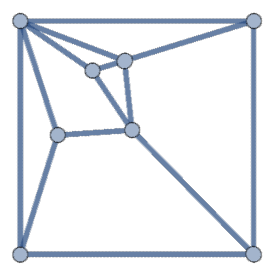}
		\includegraphics[width=\wida,clip=false]{pic34.png}
		\caption{$g_{1408.33}$ and its dual}
	\end{subfigure}
	\begin{subfigure}{\numa}
		\centering
		\includegraphics[width=\wida,clip=false]{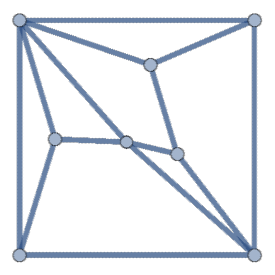}
		\includegraphics[width=\wida,clip=false]{pic11.png}
		\caption{$g_{1408.35}$ and its dual}
	\end{subfigure}
	\begin{subfigure}{\numa}
		\centering
		\includegraphics[width=\wida,clip=false]{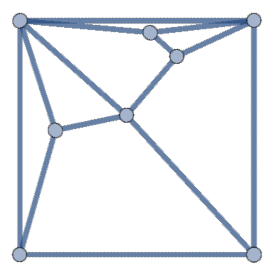}
		\includegraphics[width=\wida,clip=false]{pic21.png}
		\caption{$g_{1408.37}$ and its dual}
	\end{subfigure}
	\begin{subfigure}{\numa}
		\centering
		\includegraphics[width=\wida,clip=false]{pica31.png}
		\includegraphics[width=\wida,clip=false]{pic01.png}
		\caption{$g_{1408.39}$ and its dual}
	\end{subfigure}
	\begin{subfigure}{\numa}
		\centering
		\includegraphics[width=\wida,clip=false]{pic20.png}
		\includegraphics[width=\wida,clip=false]{pic06.png}
		\caption{$g_{1408.41}$ and its dual}
	\end{subfigure}
	\caption{The $26$ non-self-dual polyhedra with $q=14$ and $p=8$.}
	\label{fig:14.3}
\end{figure}

\clearpage
\bibliographystyle{plain}
\bibliography{bibgra}

\end{document}